\documentclass[11pt]{article}
\usepackage{amsfonts}
\usepackage{amsmath}

\newtheorem{theorem}{Theorem}[section]

\newtheorem{definition}{Definition}[section]
\newtheorem{remark}{Remark}[section]
\newtheorem{example}{Example}[section]
\numberwithin{equation}{section}
\newenvironment{proof}[1][Proof]{\noindent\textbf{#1.} }{\hfill {$\Box$}}

\begin{document}

\title{\textsc{Equivalent Definitions for Uniform Exponential Trichotomy of
Evolution Operators in Banach Spaces }}
\author{\textsc{Mihail Megan} \ \ \textsc{Codru\c{t}a Stoica}}
\date{}
\maketitle

\begin{abstract}
\noindent The aim of this paper is to give necessary and
sufficient conditions for the uniform exponential trichotomy
property of nonlinear evolution operators in Banach spaces. The
obtained results are generalizations for infinite-dimensional case
of some well-known results of Elaydi and Hajek on exponential
trichotomy of differential systems.

\vspace{3mm}

\noindent \textit{Mathematics Subject Classification:} 34D09

\vspace{3mm}

\noindent \textit{Keywords:} uniform exponential trichotomy,
evolution operator
\end{abstract}

\section{Introduction}

It is well known that in the last decades the theory of asymptotic
behaviors of evolution operators has witnessed  an explosive
development. A number of long standing open problems have recently
been solved and the theory seems to have obtained a certain degree
of maturity. There are various conditions characterizing
exponentially stable or dichotomic evolution operators on Banach or
Hilbert spaces.

In recent years, the techniques used in the investigation of the
exponential stability have been generalized for
the case of exponential dichotomy.

The concept of uniform exponential trichotomy is a natural
generalization of the classical concept of uniform exponential
dichotomy. In the study of the trichotomy, the main idea is to
obtain a decomposition of the space at every moment into three
closed subspaces: the stable subspace, the unstable subspace and
the center manifold. For the finite dimensional case some concepts
of trichotomy have been considered by Sacker and Sell in [8] and
by Elaydi and Hajek in [1], [2] and [3]. The exponential
trichotomy property in the infinite dimensional case has been
studied in [4], [5], [6] and [7].

The aim of the present paper is to give two characterizations of
the uniform exponential trichotomy of nonlinear evolution
operators on $\mathbb{R}_{+}$. We consider a concept of
exponential trichotomy which is a direct generalization of the
concept of uniform exponential dichotomy. The obtained results are
extensions for nonlinear infinite dimensional case of some well
known results of Elaydi and Hajek ([1], [2] and [3]) on
exponential trichotomy of linear differential systems.

It is important to observe that in our paper we consider a very
general concept of nonlinear evolution operators.

\section{Definitions and notations}

\noindent Let $X$ be a real or complex Banach space. The norm on
$X$ will be denoted by $\left\Vert \cdot \right\Vert $. The set of
all
mappings from $X$ into itself is denoted by $\mathfrak{F}(X)$. Let $T$ be the set of all pairs $%
(t,t_{0})$ of real numbers with $t\geq t_{0}\geq 0$.
\begin{definition}\label{eo}\rm
A mapping $E:T\rightarrow \mathfrak{F}(X)$ with the property%
\begin{equation}
E(t,s)E(s,t_{0})=E(t,t_{0}), \ \forall (t,s),(s,t_{0})\in T
\end{equation}
is called \textit{evolution operator} on $X$.
\end{definition}

\begin{example}\label{exop}\rm
If $f:\mathbb{R}\rightarrow \mathbb{R}\setminus\{0\}$ then the
mapping $E_{f}:T\rightarrow \mathfrak{F}(\mathbb{R})$ defined by
\begin{equation*}
E_{f}(t,t_{0})x=\frac{f(t)}{f(t_{0})}x
\end{equation*}
is an evolution operator on $\mathbb{R}$.
\end{example}

\begin{example}\rm
If $(S(t))_{t\geq 0}$ is a nonlinear semigroup on $X$, then the
mapping $E_{f}:T\rightarrow \mathfrak{F}(X)$ given by
$E(t,s)=S(t-s)$, where $t\geq s\geq 0$, defines an evolution
operator on $X$.
\end{example}

\begin{definition}\label{pf}\rm
An application $P:\mathbb{R}_{+}\rightarrow \mathfrak{F}(X)$ is
said to be a \textit{projection family} on $X$ if
\begin{equation}
P(t)^{2}=P(t), \ \forall t\in \mathbb{R}_{+}.
\end{equation}
\end{definition}

\begin{definition}\label{comp3}\rm
Three projection families $P_{0}$, $P_{1}$,
$P_{2}:\mathbb{R}_{+}\rightarrow
\mathfrak{F}(X)$ are said to be \textit{compatible} with the evolution operator $%
E:T\rightarrow \mathfrak{F}(X)$ if

$(c_{1})$ $P_{0}(t)+P_{1}(t)+P_{2}(t)=I$ (the identity operator)
for all $t\geq 0$;

$(c_{2})$ $P_{i}(t)P_{j}(t)=0$ for all $t\geq 0$ and all $i,j\in
\{0,1,2\},\ i\neq j$;

$(c_{3})$ $\left\Vert
P_{i}(t)x+P_{j}(t)x\right\Vert^{2}=\left\Vert
P_{i}(t)x\right\Vert^{2}+\left\Vert P_{j}(t)x\right\Vert^{2}$ for
all $t\geq 0$, all $x\in X$ and all $i,j\in \{0,1,2\},\ i\neq j$;

$(c_{4})$ $E(t,t_{0})P_{k}(t_{0})=P_{k}(t)E(t,t_{0})$ for all
$(t,t_{0})\in T$ and all $k\in \{0,1,2\}.$
\end{definition}

In what follows we will denote
\begin{equation}
E_{k}(t,t_{0})=E(t,t_{0})P_{k}(t_{0})=P_{k}(t)E(t,t_{0})
\end{equation}
for all $(t,t_{0})\in T$ and all $k\in \{0,1,2\}.$

\begin{remark}\rm
$E_{0},E_{1}$ and $E_{2}$ are evolution operators on $X.$
\end{remark}

\begin{definition}\label{uet}\rm
An evolution operator $E:T\rightarrow \mathfrak{F}(X)$ is
said to be \textit{uniformly exponentially trichotomic} if there exist some constants $%
N_{0}$, $N_{1}$, $N_{2}>1$, $\nu _{0}$, $\nu _{1}$, $\nu _{2}>0$
and three projection families $P_{0}$, $P_{1}$ and $P_{2}$
compatible with $E$ such that

$(t_{1})$ $e^{\nu _{1}(t-s)}\left\Vert E_{1}(t,t_{0})x\right\Vert
\leq N_{1}\left\Vert
E_{1}(s,t_{0})x\right\Vert $%

$(t_{2})$ $e^{\nu _{2}(t-s)}\left\Vert E_{2}(s,t_{0})x\right\Vert
\leq N_{2}\left\Vert E_{2}(t,t_{0})\right\Vert $

$(t_{3})$ $\left\Vert E_{0}(s,t_{0})x\right\Vert \leq N_{0}e^{\nu
_{0}(t-s)}\left\Vert E_{0}(t,t_{0})x\right\Vert $

$(t_{4})$ $\left\Vert E_{0}(t,t_{0})x\right\Vert \leq N_{0}e^{\nu
_{0}(t-s)}\left\Vert E_{0}(s,t_{0})x\right\Vert $

\noindent for all $t\geq s\geq t_{0}\geq 0$ and all $x\in X.$
\end{definition}

\begin{remark}\rm
In Definition \ref{uet} one can consider
\begin{equation*}
N_{0}=N_{1}=N_{2}=N \ \textrm{and} \ \nu _{1}=\nu _{2}=\nu.
\end{equation*}
Otherwise we can denote
\begin{equation*}
N=\max \left\{ N_{0},N_{1},N_{2}\right\} \ \textrm{and} \ \nu
=\min \left\{ \nu _{1},\nu _{2}\right\}.
\end{equation*}
\end{remark}

\begin{remark}\rm
For the particular case $P_{0}(t)=0$, $t\geq0$, we obtain the
uniform exponential dichotomy property. Thus, the uniform
exponential trichotomy is a natural generalization of the uniform
exponential dichotomy property.
\end{remark}

\begin{example}\rm
Let $f_{0}$, $f_{1}$, $f_{2}:\mathbb{R}\rightarrow
\mathbb{R}\setminus\{0\}$ be three functions defined by
\[
f_{0}(t)=1, \ f_{1}(t)=e^{-t}, \ f_{2}(t)=e^{t}.
\]
It is easy to observe that the evolution operators $E_{f_{0}}$,
$E_{f_{1}}$ and $E_{f_{2}}$, defined as in Example \ref{exop}, are
uniformly exponentially trichotomic.
\end{example}

\begin{example}\rm
Let us consider $X=\mathbb{R}^{3}$ with the norm
\begin{equation*}
\left\Vert(x_{1},x_{2},x_{3})\right\Vert=|x_{1}|+|x_{2}|+|x_{3}|,
\ x=(x_{1},x_{2},x_{3})\in X.
\end{equation*}
Let $\varphi:\mathbb{R}_{+}\rightarrow (0, \infty)$ be a
decreasing continuous function with the property that there exists
$\underset{t\rightarrow \infty}\lim{\varphi(t)}=l>0$.

Then the mapping $E:T\rightarrow \mathfrak{F}(X)$ defined by
\[
E(t,t_{0})x=(e^{-\int_{t_{0}}^{t}\varphi(\tau)d\tau}x_{1}, \
e^{\int_{t_{0}}^{t}\varphi(\tau)d\tau}x_{2}, \
e^{-(t-t_{0})\varphi(0)}x_{3})
\]
is an evolution operator on $X$.

Let us consider the projections defined by
\[
P_{1}(t)(x_{1},x_{2},x_{3})=(x_{1},0,0)
\]
\[
P_{2}(t)(x_{1},x_{2},x_{3})=(0,x_{2},0)
\]
\[
P_{3}(t)(x_{1},x_{2},x_{3})=(0,0,x_{3}).
\]
for all $t\geq 0$ and all $x=(x_{1},x_{2},x_{3})\in X$.

Following relations hold
\[
\left\Vert E(t,t_{0})P_{1}(t_{0})x)\right\Vert \leq
e^{-l(t-s)}\left\Vert E(s,t_{0})P_{1}(t_{0})x)\right\Vert
\]
\[
\left\Vert E(t,t_{0})P_{2}(t_{0})x)\right\Vert \geq
e^{l(t-s)}\left\Vert E(s,t_{0})P_{2}(t_{0})x)\right\Vert
\]
\[
\left\Vert E(t,t_{0})P_{3}(t_{0})x)\right\Vert \leq
e^{\varphi(0)(t-s)}\left\Vert E(s,t_{0})P_{3}(t_{0})x)\right\Vert
\]
\[
\left\Vert E(t,t_{0})P_{3}(t_{0})x)\right\Vert \geq
e^{-\varphi(0)(t-s)}\left\Vert E(s,t_{0})P_{3}(t_{0})x)\right\Vert
\]
for all $t\geq s\geq t_{0}\geq 0$ and all $x\in X$.

It follows that $E$ is uniformly exponentially trichotomic.
\end{example}

\section{The main results}
It is well known that the uniform exponential dichotomy involves
two commuting families of projections. In order to obtain a
characterization of the uniform exponential trichotomy property
using two commuting projection families we introduce the following

\begin{definition}\label{comp2}\rm
Two projection families $Q_{1}$, $Q_{2}:\mathbb{R}_{+}\rightarrow
\mathfrak{F}(X)$ are said to be \textit{compatible} with the evolution operator $%
E:T\rightarrow \mathfrak{F}(X)$ if

$(c_{1}')$ $%
Q_{1}(t)Q_{2}(t)=Q_{2}(t)Q_{1}(t)=0$

$(c_{2}')$ $\left\Vert \left[ Q_{1}(t)+Q_{2}(t)%
\right] x\right\Vert ^{2}=\left\Vert Q_{1}(t)x\right\Vert
^{2}+\left\Vert Q_{2}(t)x\right\Vert ^{2}$

$(c_{3}')$ $\left\Vert \left[ I-Q_{1}(t)\right] x\right\Vert
^{2}=\left\Vert \left[ I-Q_{1}(t)-Q_{2}(t)\right] x\right\Vert
^{2}+\left\Vert Q_{2}(t)x\right\Vert ^{2}$

$(c_{4}')$ $\left\Vert \left[ I-Q_{2}(t)\right] x\right\Vert
^{2}=\left\Vert \left[ I-Q_{1}(t)-Q_{2}(t)\right] x\right\Vert
^{2}+\left\Vert Q_{1}(t)x\right\Vert ^{2}$

$(c_{5}')$ $%
E(t,t_{0})Q_{k}(t_{0})=Q_{k}(t)E(t,t_{0}),$ $k\in \{1,2\}$

\noindent for all $t\geq 0,$ $(t,t_{0})\in T$ and all $x\in X.$
\end{definition}

\begin{remark}\rm
If $X$ is a real Hilbert space then statement $\left(
\mathit{c}_{1}^{\prime }\right)$ of Definition \ref{comp2} implies
$\left( \mathit{c}_{2}^{\prime }\right),\left(
\mathit{c}_{3}^{\prime }\right)$ and $\left(
\mathit{c}_{4}^{\prime }\right)$.
\end{remark}

The first main result of this paper is

\begin{theorem}\label{th2op}
The evolution operator $E:T\rightarrow \mathfrak{F}(X)$ is
uniformly
exponentially trichotomic if and only if there exist some constants $N>1$, $%
\nu$, $\nu _{0}>0$ and two projection families $Q_{1}$, $Q_{2}:\mathbb{R%
}_{+}\rightarrow \mathfrak{F}(X)$ compatible with $E$ such that

$(t_{1}')$ $e^{\nu (t-s)}\left\Vert
E(t,t_{0})Q_{1}(t_{0})x\right\Vert \leq N\left\Vert
E(s,t_{0})Q_{1}(t_{0})x\right\Vert $

$(t_{2}')$ $e^{\nu (t-s)}\left\Vert
E(s,t_{0})Q_{2}(t_{0})x\right\Vert \leq N\left\Vert
E(t,t_{0})Q_{2}(t_{0})x\right\Vert $

$(t_{3}')$ $\left\Vert E(s,t_{0})\left[ I-Q_{1}(t_{0})\right]
x\right\Vert \leq Ne^{\nu _{0}(t-s)}\left\Vert E(t,t_{0})\left[
I-Q_{1}(t_{0})\right] x\right\Vert $

$(t_{4}')$ $\left\Vert E(t,t_{0})\left[ I-Q_{2}(t_{0})\right]
x\right\Vert \leq Ne^{\nu _{0}(t-s)}\left\Vert E(s,t_{0})\left[
I-Q_{2}(t_{0})\right] x\right\Vert $

\noindent for all $t\geq s\geq t_{0}\geq 0$ and all $x\in X.$
\end{theorem}
\begin{proof}
\emph{Necessity}. If we denote $Q_{1}=P_{1},Q_{2}=P_{2}$ then the
conditions $(c_{1}')$, $(c_{2}')$ of Definition \ref{comp2}
respectively $(c_{5}')$ result from $(c_{2})$, $(c_{3})$
respectively $(c_{4})$ of Definition \ref{comp3}.

For $(c_{3}')$ we observe that by $(c_{3})$ we have
\begin{equation*}
\left\Vert \left[ I-Q_{1}(t)\right] x\right\Vert^{2}=\left\Vert
\left[ I-P_{1}(t)\right] x\right\Vert^{2}=\left\Vert \left[
P_{0}(t)+P_{2}(t)\right] x\right\Vert^{2}=
\end{equation*}
\begin{equation*}
=\left\Vert P_{0}(t) x\right\Vert^{2}+\left\Vert P_{2}(t)
x\right\Vert^{2}=\left\Vert \left[ I-P_{1}(t)-P_{2}(t)\right]
x\right\Vert^{2}+\left\Vert Q_{2}(t) x\right\Vert^{2}=
\end{equation*}
\begin{equation*}
=\left\Vert \left[ I-Q_{1}(t)-Q_{2}(t)\right]
x\right\Vert^{2}+\left\Vert Q_{2}(t) x\right\Vert^{2}
\end{equation*}
for all $t\geq 0$ and all $x\in X.$

Similarly for $(c_{4}')$ we have
\begin{equation*}
\left\Vert \left[ I-Q_{2}(t)\right] x\right\Vert^{2}=\left\Vert
\left[ I-P_{2}(t)\right] x\right\Vert^{2}=\left\Vert \left[
P_{0}(t)+P_{1}(t)\right] x\right\Vert^{2}=
\end{equation*}
\begin{equation*}
=\left\Vert P_{0}(t) x\right\Vert^{2}+\left\Vert P_{1}(t)
x\right\Vert^{2}=\left\Vert \left[ I-Q_{1}(t)-Q_{2}(t)\right]
x\right\Vert^{2}+\left\Vert Q_{1}(t) x\right\Vert^{2}
\end{equation*}
for all $t\geq 0$ and all $x\in X.$

Thus, the projection families $Q_{1}$ and $Q_{2}$ are compatible
with $E$.

The relations $(t_{1}')$ respectively $(t_{2}')$ result from
$(t_{1})$ respectively $(t_{2})$ of Definition \ref{uet}.

Using conditions $(c_{3})$, $(c_{4})$, $(t_{2})$ and $(t_{3})$ we
obtain
\begin{equation*}
\left\Vert E(s,t_{0})\left[ I-Q_{1}(t_{0})\right]
x\right\Vert^{2}=\left\Vert E(s,t_{0})\left[
P_{0}(t_{0})+P_{2}(t_{0})\right] x\right\Vert^{2}=
\end{equation*}
\begin{equation*}
=\left\Vert P_{0}(s)E(s,t_{0}) x\right\Vert^{2}+\left\Vert
P_{2}(s)E(s,t_{0}) x\right\Vert^{2}=\left\Vert E_{0}(s,t_{0})
x\right\Vert^{2}+\left\Vert E_{2}(s,t_{0}) x\right\Vert^{2}\leq
\end{equation*}
\begin{equation*}
\leq N^{2}e^{2\nu_{0}(t-s)}\left\Vert E_{0}(t,t_{0})
x\right\Vert^{2}+N^{2}e^{-2\nu(t-s)}\left\Vert E_{2}(t,t_{0})
x\right\Vert^{2}\leq
\end{equation*}
\begin{equation*}
\leq N^{2}e^{2\nu_{0}(t-s)}\left\Vert E(t,t_{0})\left[
P_{0}(t_{0})+P_{2}(t_{0})\right] x\right\Vert^{2}=
\end{equation*}
\begin{equation*}
=N^{2}e^{2\nu_{0}(t-s)}\left\Vert E(t,t_{0})\left[
I-Q_{1}(t_{0})\right] x\right\Vert^{2}
\end{equation*}
for all $t\geq s\geq t_{0}\geq0$ and all $x\in X$, which proves
$(t_{3}').$

The proof of $(t_{4}')$ is similar.

\emph{Sufficiency}. If we denote $P_{0}=I-Q_{1}-Q_{2},
P_{1}=Q_{1}$ and $P_{2}=Q_{2}$ then from $(c_{1}')$-$(c_{5}')$ the
statements $(c_{1})$-$(c_{4})$ are obtained immediately.

Moreover, $(t_{1}')\Leftrightarrow (t_{1})$ and
$(t_{2}')\Leftrightarrow (t_{2}).$

We observe that $P_{0}=(I-Q_{1})(I-Q_{2})$ and by $(t_{3}')$ we
obtain
\begin{equation*}
\left\Vert E_{0}(s,t_{0})x\right\Vert=\left\Vert E(s,t_{0})
P_{0}(t_{0}) x\right\Vert=\left\Vert E(s,t_{0})\left[
I-Q_{1}(t_{0})\right]\left[ I-Q_{2}(t_{0})\right] x\right\Vert\leq
\end{equation*}
\begin{equation*}
\leq Ne^{\nu_{0}(t-s)}\left\Vert E(t,t_{0})\left[
I-Q_{1}(t_{0})\right]\left[ I-Q_{2}(t_{0})\right] x\right\Vert
=Ne^{\nu_{0}(t-s)} \left\Vert E_{0}(t,t_{0}) x\right\Vert
\end{equation*}
for all $t\geq s\geq t_{0}\geq0$ and all $x\in X$. Thus $(t_{3})$
is proved.

Similarly, from $(t_{4}')$ and the remark that
$P_{0}=(I-Q_{2})(I-P_{1})$ inequality $(t_{4})$ is obtained.

It follows that $E$ is uniformly exponentially trichotomic.
\end{proof}

In order to obtain a characterization of the uniform exponential
trichotomy property using four commuting projection families, we
introduce the following

\begin{definition}\label{comp4}\rm
Four projection families $R_{1}$, $R_{2}$, $R_{3}$, $R_{4}:\mathbb{R}%
_{+}\rightarrow \mathfrak{F}(X)$ are said to be
\textit{compatible} with the evolution operator $E:T\rightarrow
\mathfrak{F}(X)$ if

$(c_{1}'')$ $%
R_{1}(t)+R_{3}(t)=R_{2}(t)+R_{4}(t)=I$

$(c_{2}'')$ $%
R_{1}(t)R_{2}(t)=R_{2}(t)R_{1}(t)=0$ and $%
R_{3}(t)R_{4}(t)=R_{4}(t)R_{3}(t)$

$(c_{3}'')$ $\left\Vert \left[ R_{1}(t)+R_{2}(t)\right]
x\right\Vert ^{2}=\left\Vert R_{1}(t)x\right\Vert ^{2}+\left\Vert
R_{2}(t)x\right\Vert ^{2}$

$(c_{4}'')$ $\left\Vert \left[ R_{1}(t)+R_{3}(t)R_{4}(t)\right]
x\right\Vert ^{2}=\left\Vert R_{1}(t)x\right\Vert ^{2}+\left\Vert
R_{3}(t)R_{4}(t)x\right\Vert ^{2}$

$(c_{5}'')$ $\left\Vert \left[ R_{2}(t)+R_{3}(t)R_{4}(t)\right]
x\right\Vert ^{2}=\left\Vert R_{2}(t)x\right\Vert ^{2}+\left\Vert
R_{3}(t)R_{4}(t)x\right\Vert ^{2}$

$(c_{6}'')$ $%
E(t,t_{0})R_{k}(t_{0})=R_{k}(t)E(t,t_{0}),$ $k\in \{1,2,3,4\}$

\noindent for all $t\geq 0,(t,t_{0})\in T$ and for all $x\in X.$
\end{definition}

\begin{remark}\rm
In the particular case when $X$ is a real Hilbert space the conditions $%
(c_{1}'')$ and $(c_{2}'')$ imply $(c_{3}''),$ $(c_{4}'')$ and
$(c_{5}'')$ in Definition \ref{comp4}.
\end{remark}

The second main result of this paper is the following

\begin{theorem}\label{th4op}
The evolution operator $E:T\rightarrow \mathfrak{F}(X)$ is
uniformly exponentially trichotomic if and only if there exist
$N>1$, $\nu$, $\nu _{0}>0$ and four projection families $R_{1}$,
$R_{2}$, $R_{3}$, $R_{4}:\mathbb{R}_{+}\rightarrow
\mathfrak{F}(X)$ compatible with $E$ such that

$(t_{1}'')$ $e^{\nu (t-s)}\left\Vert
E(t,t_{0})R_{1}(t_{0})x\right\Vert \leq N\left\Vert
E(s,t_{0})R_{1}(t_{0})x\right\Vert $

$(t_{2}'')$ $e^{\nu (t-s)}\left\Vert
E(s,t_{0})R_{2}(t_{0})x\right\Vert \leq N\left\Vert
E(t,t_{0})R_{2}(t_{0})x\right\Vert $

$(t_{3}'')$ $\left\Vert E(s,t_{0})R_{3}(t_{0})x\right\Vert \leq
Ne^{\nu _{0}(t-s)}\left\Vert E(t,t_{0})R_{3}(t_{0})x\right\Vert $

$(t_{4}'')$ $\left\Vert E(t,t_{0})R_{4}(t_{0})x\right\Vert \leq
Ne^{\nu _{0}(t-s)}\left\Vert E(s,t_{0})R_{4}(t_{0})x\right\Vert $

\noindent for all $t\geq s\geq t_{0}\geq 0$ and all $x\in X.$
\end{theorem}

\begin{proof}
\emph{Necessity}. We suppose that $E$ is uniformly exponentially
trichotomic and we denote $R_{1}=P_{1},R_{2}=P_{2},$
$R_{3}=I-P_{1},$ $R_{4}=I-P_{2}$
where $P_{0},P_{1},P_{2}$ are given by Definition \ref{uet}. Then  $%
R_{3}R_{4}=R_{4}R_{3}=P_{0}$ and the conditions
$(c_{1}'')$-$(c_{6}'')$ of Definition \ref{comp4} result
immediately from $(c_{1})$-$(c_{4})$. Thus we obtain that the
projection families $R_{1},R_{2},R_{3}$ and $R_{4}$ are compatible
with $E$.

It is obvious that $(t_{1}'')$ $%
\Leftrightarrow $ $(t_{1})$ and $(t_{2}'')$ $\Leftrightarrow $
$(t_{2}).$

In order to prove $(t_{3}'')$ we observe that
$R_{3}=I-P_{1}=P_{0}+P_{2}$ and similarly as in the proof of
$(t_{3}')$ from Theorem \ref{th2op} we obtain the desired result.
Similarly is also proved $(t_{4}'').$

\emph{Sufficiency}. If we denote $P_{1}=R_{1}$, $P_{2}=R_{2}$ and $%
P_{0}=R_{3}R_{4}$ then, from $(c_{1}'')$
and $(c_{2}'')$, we obtain%
\begin{equation*}
P_{0}+P_{1}+P_{2}=(I-R_{1})(I-R_{2})+R_{1}+R_{2}
\end{equation*}%
and hence $(c_{1})$ holds.

By $(c_{1}'')$ and $(c_{2}'')$ it follows that%
\begin{equation*}
P_{0}P_{1}=(I-R_{1})(I-R_{2})R_{1}=(I-R_{1}-R_{2})R_{1}=0
\end{equation*}%
\begin{equation*}
P_{0}P_{2}=(I-R_{1})(I-R_{2})R_{2}=(I-R_{1}-R_{2})R_{2}=0
\end{equation*}%
\begin{equation*}
P_{1}P_{2}=R_{1}R_{2}=0
\end{equation*}%
and hence the condition $(c_{2})$ is verified.

From $(c_{1}'')$-$(c_{6}'')$ it results immediately $(c_{3})$ and
thus we obtain that the projection families $P_{0},P_{1}$ and
$P_{2}$ are compatible with $E$.

It is obvious that $(t_{1}'')$ $%
\Leftrightarrow $ $(t_{1})$ and $(t_{2}'')$ $\Leftrightarrow $
$(t_{2}).$

To prove $(t_{3})$ we observe that $(t_{3}'')$ implies
\begin{equation*}
\left\Vert E_{0}(s,t_{0})x\right\Vert =\left\Vert
E(s,t_{0})P_{0}(t_{0})x\right\Vert =\left\Vert
E(s,t_{0})R_{3}(t_{0})R_{4}(t_{0})x\right\Vert \leq
\end{equation*}%
\begin{equation*}
\leq Ne^{\nu _{0}(t-s)}\left\Vert
E(t,t_{0})R_{3}(t_{0})R_{4}(t_{0})x\right\Vert =Ne^{\nu
_{0}(t-s)}\left\Vert E_{0}(t,t_{0})x\right\Vert
\end{equation*}%
for all $t\geq s\geq t_{0}\geq 0$ and all $x\in X$.

Similarly, we prove that $(t_{4}'')$ implies $(t_{4}).$

Finally, we conclude that $E$ is uniformly exponentially
trichotomic.
\end{proof}
\vspace{5mm}

{\footnotesize

\noindent {\bf Acknowledgements.} This work was completed with
financial support from the Grant CEX05-D11-23 of the Romanian
Ministry of Education and Research.

\begin{flushleft}

Mihail Megan \\

Faculty of Mathematics \\

West University of Timi\c soara \\

Romania\\

Email:\texttt{ mmegan@rectorat.uvt.ro}

\vspace{5mm}

Codru\c ta Stoica \\

Department of Mathematics \\

Aurel Vlaicu University of Arad \\

Romania\\

Email:\texttt{ stoicad@rdslink.ro}

\end{flushleft}
}

\end{document}